\documentclass{amsproc}
\usepackage{graphicx}
\usepackage{amssymb}
\usepackage{epstopdf}
\usepackage{amsmath,amsthm,amscd,amssymb}
\usepackage{latexsym}
\usepackage[colorlinks,citecolor=red,pagebackref,hypertexnames=false]{hyperref}
\usepackage{geometry}                % See geometry.pdf to learn the layout options. There are lots.
\numberwithin{equation}{section}

\theoremstyle{plain}
\newtheorem{theorem}{Theorem}[section]
\newtheorem{lemma}[theorem]{Lemma}
\newtheorem{corollary}[theorem]{Corollary}

\theoremstyle{definition}
\newtheorem{definition}[theorem]{Definition}

\newtheorem{case[theorem]}{Case}

\theoremstyle{remark}
\newtheorem{remark}[theorem]{Remark}

\numberwithin{equation}{section}

\begin{document}

\title{\parbox{14cm}{\centering{Distance graphs in vector spaces over finite fields, coloring and pseudo-randomness}}}


\author{Derrick Hart, Alex Iosevich, Doowon Koh, Steve Senger and Ignacio Uriarte-Tuero}

\begin{abstract} In this paper we systematically study various properties of the distance graph in ${\Bbb F}_q^d$, the $d$-dimensional vector space over the finite field ${\Bbb F}_q$ with $q$ elements. In the process we compute the diameter of distance graphs and show that sufficiently large subsets of $d$-dimensional vector spaces over finite fields contain every possible finite configurations.  
\end{abstract} 

\maketitle

\tableofcontents

\section{Introduction}

\vskip.125in

The distance graph in ${\Bbb F}_q^d$ is obtained by taking ${\Bbb F}_q^d$ and connecting two vertices corresponding to $x,y \in {\Bbb F}_q^d$ by an edge if $||x-y||=a$ for a fixed $a \in {\Bbb F}_q^{*}$, the multiplicative group of ${\Bbb F}_q$, where 
$$ ||x||=x_1^2+x_2^2+\dots+x_d^2.$$ 

More generally consider the set of colors $L=\{c_1, c_2, \dots, c_{q-1}\}$ corresponding to elements of ${\Bbb F}_q^{*}$. We connect two vertices corresponding to points $x,y \in {\Bbb F}_q^d$ by a $c_j$-colored edge if $||x-y||=j$. Denote the resulting family of graphs, with the implied edge and coloring sets, by $G^{\Delta}_q$ where $q$ runs over powers of odd primes. 

The main goal of this paper is a systematic study of the distance graph, including its diameter and pseudo-randomness properties. In the course of this investigation we prove sharp estimates for intersections of algebraic and non-algebraic varieties in ${\Bbb F}_q^d$ and the existence of arbitrary 
$k$ point configurations in sufficiently large subsets thereof. 

\subsection{Kaleidoscopic pseudo-randomness} We say that the family of graphs ${\{G_j\}}_{j=1}^{\infty}$ with the set of colors 
$$L=\{c_1, c_2, \dots, c_{|L|}\}$$ and the edge set ${\mathcal E}_j=\cup_{i=1}^{|L|} {\mathcal E}^i_j$, with ${\mathcal E}^i_j$ corresponding to the color $c_i$, is {\it kaleidoscopically pseudo-random} if there exist constants $C,C'>0$ such that the following conditions are satisfied: \begin{itemize} 
\item \begin{equation} \label{sizegrowth} |G_j| \to \infty \ \text{as} \ j \to \infty. \end{equation} 
\item \begin{equation} \label{uniformity} \frac{1}{C'} |{\mathcal E}^{i'}_{j}| \leq |{\mathcal E}^{i}_j| 
\leq C'|{\mathcal E}^{i'}_{j}|. 
\end{equation} 
\item $G_j$ is asymptotically complete in the sense that
\begin{equation} \label{completeness} \lim_{j \to \infty} \frac{{|G_j| \choose 2}-\sum_{i=1}^{|L|} |{\mathcal E}_j^i|}{{|G_j| \choose 2}}=0. \end{equation}
\item If $1 \leq k-1 \leq n$ and $L' \subset L$, with $|L'| \leq |L|-{k \choose 2}+n$, 
then any sub-graph $H$ of $G_j$ of size 
\begin{equation} \label{everysubgraph} \ge C{|G_j|}^{\frac{k-1}{k}} {|L|}^{\frac{n}{k}}, \end{equation} contains every possible sub-graph with $k$ vertices and $n$ edges with an arbitrary edge color distribution from $L'$. \end{itemize}

See, for example, a survey by Krivelevich and Sudakov (\cite{KS07}) for related notions of pseudo-random graphs, examples and applications. The first result of this paper is the following. 
\begin{theorem} \label{maingraph} The above defined family of graphs $\{G^{\Delta}_q\}$ is kaleidoscopically pseudo-random. \end{theorem} 

The proof shows that the constant $C'$ in the definition of kaleidoscopic pseudo-randomness may be taken to be $(1+o(1))$ in this context if the dimension $d$ is not two or 
the zero distance is excluded. The constant $C$ that the proof yields is exponential in the number of vertices. 

We actually prove a little more as the arguments below indicate. We shall see that under the set of hypotheses corresponding to kaleidoscopic pseudo-randomness, every finite geometric configuration in ${\Bbb F}_q^d$ is realized. See \cite{V07} and \cite{Vu07} where related questions are studied using graph theoretic methods. 

The first item in the definition of weak pseudo-randomness above (\ref{sizegrowth}) is automatic as the size of $G_j$ is $q^d$, by construction. The second and third items, (\ref{uniformity}) and (\ref{completeness}), respectively, are easy special cases of the following calculation.  While it is implicit in (\cite{IR07}), we give the proof at the end of the paper for the sake of reader's convenience. 
\begin{lemma} \label{colorsize} For any $t \in {\Bbb F}_q$, 
$$ |\{(x,y) \in {\Bbb F}_q^d \times {\Bbb F}_q^d: ||x-y||=t\}|=\left\{\begin{array}{ll} (2+o(1))\,q^{2d-1} \quad\mbox{if} \quad d=2, t=0\\
                                                                                        (1+o(1))\,q^{2d-1} \quad\mbox{otherwise} \end{array}\right.$$ 
where $o(1)$ means that the quantity goes to $0$ as $q \to \infty$. 
\end{lemma} 

We now ready to address the meat of our definition of weak pseudo-randomness, which is the fourth item (\ref{everysubgraph}). 
 
\begin{definition} Given $L' \subset {\Bbb F}_q^{*}$ such that 
$$|L'| \leq q-1-{k \choose 2}+|J|,$$ and 
$$J \subset {\{1,2, \dots, k\}}^2 \ \backslash \ \{(i,i): 1 \leq i \leq k\},$$ a $k$-point $J$-configuration in 
$E$ is a set of $k$ points $\{x^1, x^2, \dots, x^k\}$ such that 
$$||x^i-x^j||=a_{ij} \in L'$$ for all $(i,j) \in J$. Denote the set of all $k$ point $J$-configurations by ${\mathtt T}_k^J(E)$. \end{definition} The item (\ref{everysubgraph}) follows from the following geometric estimate. 
\begin{theorem} \label{main} Let $E \subset {\Bbb F}_q^d$, $d \ge 2$. Suppose that $1 \leq k-1 \leq n \leq d$ and 
\begin{equation} \label{sizeconfig} |E| \ge Cq^{d \left( \frac{k-1}{k}\right)} q^{\frac{n}{k}} \end{equation} with a sufficiently large constant $C>0$. Then for any 
$$J \subset {\{1,2, \dots, k\}}^2 \ \backslash \ \{(i,i): 1 \leq i \leq k\}$$ with $|J|=n$, we have 
$$ |{\mathtt T}_k^J(E)|=(1+o(1)){|E|}^{k}q^{-n}.$$ 
\end{theorem} 

Our proof uses geometric and character sum machinery similar to the one used in \cite{IR07} and \cite{HI07}. In the former paper, Theorem \ref{main} is proved in the case $k=2$ and $n=1$, and in the latter article Theorem \ref{main} is demonstrated in the case of general $k$ and $n={k \choose 2}$. Thus Theorem \ref{main} and, consequently, Theorem \ref{maingraph} may be viewed as filling the gap between these results. 

\vskip.125in 

\subsection{Diameter of the distance graph and related objects} 

\vskip.125in 

Let the distance graph $G_q^{\Delta}$, equipped with the coloring set $L$ be as above. Given a fixed color $c$ in $L$, we define the diameter of $G_q^{\Delta}$ as follows. Given vertices $x$, $y$ in $G_q^{\Delta}$, define a {\it path} of length $k$ from $x$ to $y$ to be a sequences $\{x^1, \dots, x^{k+1}\}$, where $x^j$s are distinct, $x^1=x$, $x^{k+1}=y$, each $x^j$ is a vertex in $G_q^{\Delta}$ and $x^i$ is connected to $x^{i+1}$ by a $c$-colored edge for every $1 \leq i \leq k$. We say that a path from $x$ to $y$ is optimal if it is a path and its length is as small as possible. Define the {\it diameter} of $G_q^{\Delta}$, with respect to the color $c$, to be the largest length of the optimal path between any two vertices in $G_q^{\Delta}$. 

Our first result in this direction is actually about a more general families of graphs. Let $U \subset {\Bbb F}_q^d$. We say that $U$ is Salem if there exists a uniform constant $C>0$ such that 
$$ |\widehat{U}(\xi)| \leq Cq^{-d} {|U|}^{\frac{1}{2}},$$ where the Fourier transform with respect to a non-trivial additive character $\chi$ is defined and briefly reviewed in (\ref{ftdef}) and the lines that follow. We shall also see below (Lemma \ref{kloosterman}) that the sphere 
\begin{equation}\label{defsphere} S_t=\{x \in {\Bbb F}_q^d: x_1^2+\dots+x_d^2=t\}\end{equation} is a Salem set. 

Define $G^U_q$ to the graph with vertices in ${\Bbb F}_q^d$ and two vertices, corresponding to $x, y \in {\Bbb F}_q^d$ connected by an edge if $x-y \in U$. We do not attach a coloring scheme in this context. 
\begin{theorem} \label{generaldiameter} Suppose that $U$ is Salem and $|U| \ge Cq^{\frac{2d}{3}}$ with a sufficiently large constant $C>0$. Then the diameter of $G^U_q$ is $\leq 3$. 
\end{theorem} 

\begin{corollary} \label{sphereweak} Given any fixed color $c \in L$, the graph $G_q^{\Delta}$ has diameter $\leq 3$ if $d \ge 4$. 
\end{corollary} 

The fact that the sphere is Salem, mentioned above, is proved in Lemma \ref{kloosterman} below as is the fact that $|S| \approx q^{d-1}$. 
It follows that the diameter of $G_q^{\Delta}$ is $\leq 3$ provided that 
$|S| \ge Cq^{\frac{2d}{3}}$ with a sufficiently large constant $C>0$. Since $|S| \approx q^{d-1}$, this holds if $d \ge 4$, which completes the proof of the corollary. 
We can do a bit better, however. 
\begin{theorem} \label{spherediametersharp} 1) If $d\geq 4$ then the diameter of $G_q^{\Delta}$ is two for all $q \ge 3.$ \\
                                            2) If $d=2$ then the diameter of $ G_q^{\Delta}$ is never two for all $q\geq 5.$
                                               Moreover, the diameter of $G_q^{\Delta}$ is three if $q\not=3,5,9,13.$ \\
                                            3) If $d=3$ then the diameter of $G_q^{\Delta}$ is two or three.\end{theorem}
\vskip.125in 

\section{Pseudo-arithmetic progressions} 

\vskip.125in 

Consider a sequence of $k$ points $P_1, P_2, \dots, P_k$ in ${\Bbb F}_q^d$ such that 
$$||P_j-P_i||={(j-i)}^2 \ \text{for} \ 1 \leq i \leq j \leq k.$$

We call such an ordered sequence of vectors {\it pseudo-arithmetic}. The following is a simple consequence of Theorem \ref{main}. 
\begin{corollary} \label{pizdatayaprogressiya} Suppose that $E \subset {\Bbb F}_q^d$ such that $|E| \ge Cq^{\frac{k-1}{k}d} q^{\frac{k-1}{2}}$. Then $E$ contains a pseudo-arithmetic progression of length $k$. \end{corollary} 

In fact, Theorem \ref{main} implies that $E$ contains $\approx {|E|}^kq^{-{k \choose 2}}$ arithmetic-like progressions. It would be wonderful if these were actual arithmetic progressions. In fact, suppose it were true that every arithmetic-like progression is an actual arithmetic progression in at least one coordinate. We could then take $E=A \times A \times \dots \times A$ and conclude that if $|A| \ge Cq^{\frac{k-1}{k}} q^{\frac{k-1}{2d}}$, then $A$ contains an arithmetic progression of length $k$, thus giving us a rather attractive version of Szemeredi's theorem in finite fields. The reality is very different, however. It is easy enough to construct examples of sequences which are arithmetic-like but not actually arithmetic. Let $z \in {\Bbb F}_q^{d-1}$ such that $||z||=z_1^2+\dots+z_{d-1}^2=0$. Let $P_j=(j,z) \in {\Bbb F}_q^d$. It is not hard to see that $||P_j-P_i||={(j-i)}^2+||z||={(j-i)}^2,$ so the sequence is arithmetic-like, but it is certainly not in general an arithmetic progression. 

What is somewhat more difficult is to construct examples of arithmetic-like sequences that are not arithmetic progressions in any coordinate. One way is to take one of the arithmetic-like progressions described in the previous paragraph and rotate it. For example, we may start out with the sequence 
$$ (0,0,0) \ (1,1,i) \ (2,0,0),$$ where $i=\sqrt{-1}$ and rotate it by an orthogonal matrix 
$$ \left( \begin{array}{ccc} 
t & -t & 0 \\
t & t & 0 \\
0 & 0 & 1 \end{array} \right).$$

In order to have the determinant of this matrix equal to $1$ we must have $t^2=1/2$. This equation has a solution in some fields and not others. Recall that we are also using $i=\sqrt{-1}$, an object which exists in some fields and not others. The simplest field where both objects exist is ${\Bbb Z}_{17}$. In this field we may take $t=3$ and $i=4$. We thus obtain the sequence 
$$ (0,0,0) \ (0,6,4) \ (6,6,0).$$ 

Observe that this sequence is not arithmetic in any coordinate. 

\vskip.125in 

\section{Proof of the "kaleidoscopic" result (Theorem \ref{main})}. 

\vskip.125in 

Let $T_k^J$ denote the set of  $k$-point $J$-configurations in $E$ and let 
$T_k^J(x^1, \dots, x^k)$ denote its characteristic function. Assume, inductively, that for every $J' \subset J$, 
\begin{equation} \label{inductivehypothesis} |T_{k-1}^{J'}|=(1+o(1)){|E|}^{k-1}q^{-|J'|} \end{equation} if 
$$ |E| \ge Cq^{d \left( \frac{k-2}{k-1} \right)} q^{\frac{|J'|}{k-1}}.$$ 

The initialization step is the following. Observe that 
$$ |T_1^J|=|T_1^{\emptyset}|=|E|=|E|q^{-0},$$ and this needs to hold if 
$$ |E| \ge C{|E|}^{d \left( \frac{k-1}{k} \right)} q^{\frac{|J|}{k}}=C.$$ 

\subsection{The induction step:} We have, without loss of generality,  
\begin{equation} \label{setup} |T_k^J|=\sum T_{k-1}^{J'}(x^1, \dots, x^{k-1}) E(x^k) \Pi_{j=1}^l S(x^j-x^k) \Pi_{i=l+1}^{k-1} \sum_{a_i \not=0} S_{a_i}(x^i-x^k) \end{equation} for some $1 \leq l \leq k-1$, depending on the degree of the vertex corresponding to $x^k$, where 
$$ S_t=\{x \in {\Bbb F}_q^d: x_1^2+x_2^2+\dots+x_d^2=t\},$$ and $S \equiv S_1$. Technically, we should replace $ \Pi_{j=1}^l S(x^j-x^k)$ by $\Pi_{j=1}^l S_{a_j}(x^j-x^k)$ for an arbitrary set of $a_j$s, but this does not change the proof any and only complicates the notation. 

Recall that given a function $f: {\Bbb F}_q^m \to {\Bbb F}_q$, the Fourier transform with respect to a non-trivial additive character $\chi$ on ${\Bbb F}_q$ is given by the relation 
\begin{equation} \label{ftdef} \widehat{f}(\xi)=q^{-m} \sum_{x \in {\Bbb F}_q^m} \chi(-x \cdot \xi). \end{equation} Also recall that 
\begin{equation} \label{inversion} f(x)=\sum_{\xi \in {\Bbb F}_q^m} \chi(x \cdot m) 
\widehat{f}(m) \end{equation} and 
\begin{equation} \label{plancherel} \sum_{\xi \in {\Bbb F}_q^m} {|\widehat{f}(\xi)|}^2=q^{-m} \sum_{x \in {\Bbb F}_q^m} {|f(x)|}^2. \end{equation} 

We shall also need the following estimates based on classical Gauss and Kloosterman sum bounds. See, for example, Lemma 2.2 and its proof in \cite{IR07} for the first and the third estimates below. 
\begin{lemma} \label{kloosterman} With the notation above, for any $t \not=0$, $\xi \not=(0, \dots, 0)$ and $q$ sufficiently large, 
\begin{equation} \label{decay} |\widehat{S}_t(\xi)| \leq 2q^{-\frac{d+1}{2}}. \end{equation}

Moreover, for any $a \not=0$, 
\begin{equation} \label{averagedecay} \left|\sum_{t \not=a} \widehat{S}_t(\xi)\right| \leq 
(2+o(1))q^{-\frac{d+1}{2}} \end{equation} and 
\begin{equation} \label{zeroterm}  \widehat{S}_t(0, \dots, 0)=q^{-d}|S_t|=(1+o(1))q^{-1},\end{equation} where $o(1)$ means that the quantity goes to $0$ as $q \to \infty$. \end{lemma} 

Using (\ref{inversion}) and the definition of the Fourier transform, we see from (\ref{setup}) that 
$$ |T_k^J|=q^{kd} \sum_{\xi^1, \dots, \xi^{k-1}: \xi^s \in {\Bbb F}_q^d} 
\widehat{T_{k-1}^{J'}}(\xi^1, \dots, \xi^{k-1}) \widehat{E} \left( \sum_{u=1}^{k-1} \xi^u \right) \Pi_{j=1}^l \widehat{S}(\xi^j) \Pi_{i=l+1}^{k-1} \sum_{a_i \not=0} \widehat{S}_{a_i}(\xi^i)$$
$$=Main+Remainder,$$ where Main is the term corresponding to taking $\xi^s=(0, \dots, 0)$ for every $1 \leq s \leq k-1$. It follows by Lemma \ref{kloosterman} that 
$$ Main=(1+o(1))|T_{k-1}^{J'}||E|q^{-l}.$$ 

The Remainder is the sum of terms of the form $R_{U,V}$, where 
$$ U=\{j \in \{1,2, \dots, l\}: \xi^j \not=(0, \dots, 0)\},$$ and 
$$ V=\{j \in \{l+1. \dots, k-1\}: \xi^j \not=(0, \dots, 0)\}.$$ 

We first analyze the term where compliments of $U$ and $V$ are empty sets. We get 
$$ R_{U,V}=q^{kd} \sum_{\xi^1, \dots, \xi^{k-1}: \xi^s \in {\Bbb F}_q^d; \xi^s \not=(0, \dots, 0)} 
\widehat{T_{k-1}^{J'}}(\xi^1, \dots, \xi^{k-1}) \widehat{E} \left( \sum_{u=1}^{k-1} \xi^u \right) \Pi_{j=1}^l \widehat{S}(\xi^j) \Pi_{i=l+1}^{k-1} \sum_{a_i \not=1} \widehat{S}_{a_i}(\xi^i).$$

Applying Lemma \ref{kloosterman} to the Fourier transforms of spheres and applying Cauchy-Schwartz, in the variables $\xi^1, \dots, \xi^{k-1}$, followed by (\ref{plancherel}) to the first two terms in the sum, we see that 
$$ R_{U,V}=O \left( q^{kd} \cdot {|T_{k-1}^{J'}|}^{\frac{1}{2}} \cdot {|E|}^{\frac{1}{2}} \cdot q^{-\frac{d}{2}} 
\cdot q^{\frac{d(k-2)}{2}} \cdot q^{-\frac{d+1}{2}l} \cdot q^{-\frac{d+1}{2}(k-1-l)} \right)$$
$$=O \left( {|T_{k-1}^{J'}|}^{\frac{1}{2}} \cdot {|E|}^{\frac{1}{2}} \cdot q^{\frac{d(k-1)}{2}} q^{-\frac{l}{2}} q^{-\frac{(k-1-l)}{2}} \right),$$ where $X=O(Y)$ means that there exists $C>0$, independent of $q$, such that $X \leq CY$. 

Applying the inductive hypothesis (\ref{inductivehypothesis}) and noting that $l$ may be as large as 
$k-1$, we see that 
$$ R_{U,V} \leq \frac{1}{2} \cdot Main$$ if 
$$ |E| \ge Cq^{d \left( \frac{k-1}{k} \right)} \cdot q^{\frac{|J|}{k}},$$ with $C$ sufficiently large, as desired. 

To estimate the general $R_{U,V}$, we need the following simple observation that is proved by a direct calculation. Let $|U|+|V|=m$ and define 
\begin{equation} \label{zeroes} \widehat{f}(\mu^1, \dots, \mu^m)=q^{d(k-1)} q^{-md} 
\widehat{T}_{k-1}^{J'}(Z_{U,V}(\mu^1), \dots, Z_{U,V}(\mu^{k-1}))), \end{equation} where 
$$ Z_{U,V}: {\Bbb F}^{d} \to {\Bbb F}^{d}$$ with $Z_{U,V}(\xi^j)=\xi^j$ if $j \in U \cup V$ and $(0, \dots, 0)$ otherwise. Then 
$$ \sum_{y^1, \dots, y^m} f^2(y^1, \dots, y^m) \leq \max_{y^1, \dots, y^m} f(y^1, \dots, y^m) \cdot \sum_{y^1, \dots, y^m} f(y^1, \dots, y^m)$$
$$ \leq \min \{|E|, (1+o(1))q^{d-1} \} \cdot |T^{J'}_{k-1}|.$$ 

Applying (\ref{zeroes}) one can check that in the regime $|E| \ge Cq^{d \left( \frac{k-1}{k} \right)} q^{\frac{n}{k}}$ the remaining $R_{U,V}$s are smaller than the error term we already estimated. This completes the proof. 

Technically speaking we must still show that if $|T_k^{J_1}|$ satisfies the conjectured estimate for every $J_1$ with $|J_1|=n_1$, then so does $T_k^{J_2}$ with $|J_2|=n_2>n_1$. However, this is apparent from the proof above.

\vskip.125in 

\section{Results based on Classical Gauss sums } 

\vskip.125in 
In this section, we collect the well-known facts which follow by estimates of Gauss sums.
Such facts shall be used in the next sections.
Let $\chi$ be a non-trivial additive character of ${\mathbb F}_q$ and $\psi$ a multiplicative character of 
${\mathbb F}_q$ of order two, that is, $ \psi(ab)=\psi(a)\psi(b)$ and  $\psi^2(a)=1$ for all $a,b \in {\mathbb F}_q^*$ but $\psi\not\equiv 1.$
For each $a\in {\mathbb{F}_q}$, the Gauss sum $G_a(\psi, \chi)$ is defined by 
$$ G_a(\psi, \chi) = \sum_{s\in {\mathbb F}_q^*} \psi(s) \chi(as).$$ The magnitude of the Gauss sum is given by the relation
$$ |G_a(\psi, \chi)| = \left\{\begin{array}{ll} q^{\frac{1}{2}} \quad &\mbox{if} \quad a\ne 0\\
                                                  0 \quad & \mbox{if} \quad a=0. \end{array}\right.$$
\begin{remark}
Here, and throughout this paper, we denote by $\chi$ and $\psi$ the canonical additive character  and the quadratic character of 
${\mathbb F}_q^*$ respectively. Recall that if $\psi$ is the quadratic character of ${\mathbb F}_q^*$ then 
 $\psi(s)=1$ if $s$ is a square number in ${\mathbb F}_q^*$ and $\psi(s)=-1$ otherwise.
\end{remark}
The following theorem provided us of the explicit formula of the Gauss sum $G_1(\psi, \chi).$ 
For the nice proof, see \cite{LN97}.
\begin{theorem}\label{ExplicitGauss}
Let ${\mathbb F}_q$ be a finite field with $ q= p^l$, where $p$ is an odd prime and $l \in {\mathbb N}.$
Then we have
$$G_1(\psi, \chi)= \left\{\begin{array}{ll}  {(-1)}^{l-1} q^{\frac{1}{2}} \quad &\mbox{if} \quad p =1 \,\,( mod\,4) \\
                    {(-1)}^{l-1} i^l q^{\frac{1}{2}} \quad &\mbox{if} \quad p=3 \,\,( mod\,4).\end{array}\right. $$
 \end{theorem}

In particular, we have
\begin{equation}\label{square} \sum_{s \in {\mathbb F}_q} \chi (a s^2) = {\psi(a)} G_1(\psi, \chi) \quad \mbox{for any} \quad a \ne 0,\end{equation}
because the quadratic character $\psi$ is the multiplicative character of ${\mathbb F}_q^* $ of order two. For the nice proof for this equality and
the magnitude of Gauss sums, see \cite{LN97} or \cite{IK04}. As the direct application of the equality in (\ref{square}), we have the following estimate.
\begin{lemma}\label{complete}
For $\beta \in {\mathbb F}_q^k$ and $t\ne 0$, we have
$$ \sum_{\alpha \in {\mathbb F}_q^k} \chi( t \alpha \cdot \alpha + \beta \cdot \alpha ) 
= \chi\left( \frac{\|\beta\|}{-4t}\right) \eta^k(t) \left(G_1(\eta, \chi)\right)^k,$$
where, here and throughout the paper, $\|\beta\| = \beta \cdot \beta.$
\end{lemma}
\begin{proof}
It follows that
$$\sum_{\alpha \in {\mathbb F}_q^k} \chi( t \alpha \cdot \alpha + \beta \cdot \alpha ) 
              =\prod_{j=1}^{k} \sum_{\alpha_j\in {\mathbb F}_q} \chi( t\alpha_j^2 + \beta_j\alpha_j).$$
Completing the square in $\alpha_j$-variables, changing of variables, $ \alpha_j+\frac{\beta_j}{2t} \to \alpha_j$,
and using the inequality in (\ref{square}), the proof immediately follows.\end{proof}

Due to the explicit formula for the Gauss sum $G(\psi,\chi)$, we can count the number of the elements 
in spheres $S_t \subset {\mathbb F}_q^d$ defined as before. The following theorem enables us to see the exact number of 
the elements of spheres $S_t$ which depends on the radius $t$, dimensions, and the size of the underlining finite field ${\mathbb F}_q.$

\begin{theorem}\label{explicit}
Let $S_t \subset {\mathbb F}_q^d$ be the sphere defined as in (\ref{defsphere}). For each $t \not=0,$ we have
$$ |S_t| = \left\{\begin{array}{ll} q^{d-1}- q^{\frac{d-2}{2}} \psi\left((-1)^{\frac{d}{2}}\right) \quad &\mbox{if} \quad d \quad \mbox{is even}\\
                                    q^{d-1}+ q^{\frac{d-1}{2}} \psi\left((-1)^{\frac{d-1}{2}} t\right) \quad & \mbox{if} \quad d \quad \mbox{is odd} \end{array}\right.$$
                               
\end{theorem}

\begin{proof}
For each $t\not=0,$ we have
\begin{align*} |S_t| =& \sum_{\|x\| =t} 1\\
                     =& q^{d-1} + q^{-1} \sum_{s\not=0}\sum_{x} \chi\left(s(\|x\|-t)\right).\end{align*}
Using Lemma \ref{complete}, we see that
$$ |S_t|= q^{d-1} +q^{-1} (G_1(\psi,\chi))^d \sum_{s \not=0} \psi^d(s) \chi(-st)$$

Case 1: $d$ is even. Then $ \psi^d  \equiv 1,$ because the quadratic character $\psi$ is a multiplicative character of ${\mathbb F}_q^*$ of order two.
Thus we see that
\begin{align}\label{case1} |S_t| &= q^{d-1}+q^{-1}(G_1(\psi,\chi))^d \sum_{s \not=0} \chi(-st)\nonumber\\
                    &=q^{d-1}-q^{-1}(G_1(\psi,\chi))^d .\end{align}

Case 2: $d$ is odd. Then $\psi^d = \psi,$ because the order of $\psi$ is also two. It follows that
\begin{align}\label{case2} 
|S_t|=& q^{d-1} +q^{-1} (G_1(\psi,\chi))^d \sum_{s \not=0} \psi(s) \chi(-st)\nonumber\\
     =& q^{d-1} +q^{-1} (G_1(\psi,\chi))^{d+1} \psi(-t). \end{align}
Together with (\ref{case1}) and (\ref{case2}), it suffices to show that
$$\left\{\begin{array}{ll} (G_1(\psi,\chi))^d=q^{\frac{d}{2}} \psi\left((-1)^{\frac{d}{2}}\right) \quad \mbox{if}\quad d \quad \mbox{is even}\\
(G_1(\psi,\chi))^{d+1}=q^{\frac{d+1}{2}} \psi\left((-1)^{\frac{d+1}{2}}\right) \quad \mbox{if}\quad d \quad \mbox{is odd}\end{array}\right.$$
However this follows by Theorem \ref{ExplicitGauss} and the well-known fact that
$ \psi(-1)=1$ if $ q \equiv 1 \,\, (mod \,4)$ and $\psi(-1)=-1$ if $q \equiv 3 \,\, (mod\, 4).$ Thus the proof is complete.
\end{proof}

\vskip.125in 

\section{Proof of the uniformity of color distribution (Lemma \ref{colorsize})} 

\vskip.125in 

We have 
$$ |\{(x,y) \in {\Bbb F}_q^d \times {\Bbb F}_q^d: ||x-y||=t\}|$$
$$=\sum_{x,y} S_t(x-y)=\sum_{x,y} \sum_m \chi((x-y) \cdot m) \widehat{S}_t(m)$$
$$=q^{2d} \widehat{S}_t(0, \dots, 0)=q^d |S_t|.$$ 

If $t\not=0,$ Theorem \ref{explicit} yields that 
$$ |S_t|=(1+o(1))q^{d-1}.$$ 
On the other hand, we have
$$ |S_0|=\sum_x S_0(x)=q^{-1} \sum_x \sum_s \chi(s||x||)$$
$$=q^{d-1}+q^{-1}\sum_x \sum_{s \not=0} \chi(s||x||),$$ 
$$= q^{d-1}+q^{-1} \left(G_1(\psi, \chi)\right)^d  \sum_{s \not=0} \psi^d(s)$$ 
Since $\psi^d\equiv 1$ if $d$ is even and $ \psi^d=\psi$ if $d$ is odd, we obtain that
$$|S_0|= \left\{\begin{array}{ll}q^{d-1} +q^{-1} (q-1)\left(G_1(\psi, \chi)\right)^d \quad &\mbox{if} \quad d \quad \mbox{is even}\\
                                 q^{d-1}   &\mbox{if} \quad d \quad \mbox{is odd}\end{array}\right.$$
Thus the proof immediately follows by this, because the magnitude of the Gauss sum $G_1(\psi,\chi)$ is exactly $q^{\frac{1}{2}}.$

\vskip.125in 

\section{Proof of the Fourier decay estimates (Lemma \ref{kloosterman})} 

\vskip.125in 

We use a part of the argument above. For each $m \not=(0,\ldots,0)$, we have 
\begin{align}\label{startaverage} \widehat{S}_t(m)=&q^{-d} q^{-1} \sum_s \sum_x \chi(-x \cdot m) \chi(s(||x||-t))\nonumber \\
=&q^{-d} q^{-1} \sum_{s \not=0} \sum_x \chi(-x \cdot m) \chi(s(||x||-t))\nonumber \\
=&q^{-d-1} (G_1(\psi,\chi))^d\sum_{s \not=0} \chi \left( \frac{||m||}{-4s}-st \right) \psi^d(s),\end{align} 
and the estimate (\ref{decay}) follows from the following classical estimate due to Andre Weyl (\cite{We48}). 
\begin{theorem}\label{generalKloosterman} Let 
$$ K(a)=\sum_{s \not=0} \chi(as^{-1}+s) \phi(s),$$ where $\phi$ is a multiplicative character on ${\Bbb F}_q^{*}$. Then for any $a \in {\Bbb F}_q$, 
$$ |K(a)| \leq 2\sqrt{q}.$$ 
\end{theorem} 

We now turn our attention to (\ref{averagedecay}). We may assume that $a=1$ without loss of generality. With (\ref{startaverage}) as the starting point, 
we sum this expression in $t \not=1$ and obtain 
$$ \sum_{t\not=1} \widehat{S_t}(m)=-q^{-d-1} (G_1(\psi,\chi))^d\sum_{s \not=0} \chi \left( \frac{||m||}{-4s} \right) \psi^d(s) \chi(-s).$$

We see that this expression is 
$$ \leq 2q^{-\frac{d+1}{2}},$$ because the magnitude of the Gauss sum $G_1(\psi,\chi)$ is $q^{\frac{1}{2}}$ and the sum in $s\not =0$ is just the generalized Kloosterman sum
in Theorem \ref{generalKloosterman}.The proof of (\ref{zeroterm}) follows easily from the proof of Lemma \ref{colorsize}. This completes the proof of Lemma \ref{kloosterman}. 

\vskip.125in 

\section{Proof of Theorem \ref{generaldiameter}} 

\vskip.125in 

We shall deduce Theorem \ref{generaldiameter} from the following estimate. 
\begin{theorem} \label{twodistance} $U, E,F \subset {\Bbb F}_q^d$ such that $U$ is Salem and 
$$ |E||F| \ge C\frac{q^{2d}}{|U|}$$ with a sufficiently large constant $C>0$. Then 
$$ \nu_U=\{(x,y) \in E \times F: x-y \in U\}>0.$$ 
\end{theorem} 

Taking Theorem \ref{twodistance} for granted, for a moment, Theorem \ref{generaldiameter} follows instantly. Indeed, take $x,y$ with $x \not=y$ in ${\Bbb F}_q^d$. Let $E=U+x$ and $F=U+y$. It follows that $|E|=|F|=|U|$, so $|E||F|={|U|}^2$. We conclude from Theorem \ref{twodistance} that if $|U| \ge Cq^{\frac{2d}{3}}$ with a sufficiently large constant $C>0$, then there exists $x' \in U+x$ and $y' \in U+y$ such that $x'-y' \in U$. This implies that the diameter of $G_q^U$ is at most three as desired. 

To prove Theorem \ref{twodistance}, observe that 
$$ \nu_U=\sum_{x,y} E(x)F(y)U(x-y)$$
$$=\sum_{x,y} \sum_m \widehat{U}(m)\chi((x-y) \cdot m) E(x)F(y)$$
$$=q^{2d} \sum_m \overline{\widehat{E}(m)} \widehat{F}(m) \widehat{U}(m)$$
$$=|E||F||U|q^{-d}+q^{2d} \sum_{m \not=(0, \dots, 0)} 
\overline{\widehat{E}(m)} \widehat{F}(m) \widehat{U}(m)=I+II.$$ 

By assumption, 
$$ |II| \leq q^{2d} \cdot q^{-d} {|U|}^{\frac{1}{2}} \cdot \sum_{m \not=(0, \dots, 0)} |\overline{\widehat{E}(m)}| |\widehat{F}(m)|$$
$$ \leq Cq^d {|U|}^{\frac{1}{2}} {\left( \sum {|\widehat{E}(m)|}^2 \right)}^{\frac{1}{2}} \cdot {\left( \sum {|\widehat{F}(m)|}^2 \right)}^{\frac{1}{2}}$$
$$=C{|U|}^{\frac{1}{2}} {|E|}^{\frac{1}{2}} {|F|}^{\frac{1}{2}}$$ by (\ref{plancherel}). Comparing $I$ and $II$ we complete the proof of Theorem \ref{twodistance}. 

\vskip.125in 

\section{Proof of Theorem \ref{spherediametersharp}} 
In this section, we provide the proof of Theorem \ref{spherediametersharp}.
The proof of the first part in Theorem \ref{spherediametersharp}  is given in the following subsection \ref{diametertwo}. 
For the proof of the second and third part in Theorem \ref{spherediametersharp}, we first show in Subsection \ref{nevertwo} that 
the diameter of $G_q^{\Delta}$ in two dimension is never two if $q\not=3$ and then we complete in Subsection \ref{lastproof} 
the proof of the second and third part of Theorem \ref{spherediametersharp}.
\vskip.125in 

\subsection{The Proof of the first part of Theorem \ref{spherediametersharp}}\label{diametertwo}
We first prove that in dimensions four and higher, the diameter is two  though we will actually prove a much stronger statement. 
It suffices to show that if the dimensions $d\geq 4$, then two different spheres in ${\mathbb F}_q^d$ with same radius $t\not=0$  always intersect.
The proof is based on the following lemma.
\begin{lemma}\label{intersection2}
For each $x\not=(0,\ldots,0)$, and $t\not=0$, we have
$$ |S_t \cap (S_t +x)| =q^{-d} |S_t|^2-q^{-1}+ q^{-2} \left(G_1(\psi, \chi)\right)^{d+1} \psi(-1) \sum_{\substack{r\not=0, 1\\ :t(1-r)^2+r\|x\|\not=0}}
\psi(t(1-r)^2+r\|x\|)$$
if $d$ is odd.
On the other hand, if $d$ is even then we have
$$ |S_t \cap (S_t +x)| =q^{-d} |S_t|^2-q^{-2} - q^{-2}(q-2) \left(G_1(\psi, \chi)\right)^d +q^{-1}\left(G_1(\psi, \chi)\right)^d  
\sum_{\substack{r\not=0, 1\\: t(1-r)^2+r\|x\|=0}}1.$$
\end{lemma} 
By assuming Lemma \ref{intersection2} for a moment, we shall prove that the diameter in dimensions four and higher is two.
It suffices to show that if $x \not= (0, \ldots, 0)$ then $|S_t \cap (S_t +x)| >0$ for all $t\not=0$ and $d\geq 4.$\\

Case 1: Suppose that $d\geq 5$ is odd. Then  we see from Theorem \ref{explicit} that 
\begin{equation} \label{size1}
|S_t| \geq q^{d-1} - q^{\frac{d-1}{2}}.
\end{equation} 
On the other hand, it is clear that 
\begin{equation}\label{size2}
0\leq  \sum_{r\not=0, 1: t(1-r)^2+r\|x\|\not=0}1\leq q-2.
\end{equation}
Using the first part of Lemma \ref{intersection2} together with (\ref{size1}), (\ref{size2}), and the magnitude of the Gauss sum $G_1(\psi, \chi)$, 
we see that 
$$ |S_t \cap (S_t +x) | \geq q^{-d} \left(q^{d-1} - q^{\frac{d-1}{2}}\right)^2 -q^{-1} - q^{-2} q ^{\frac{d+1}{2}} (q-2) = q^{d-2}-q^{\frac{d-1}{2}}, $$
which is greater than zero if $d \geq 5$  as wanted.\\

Case 2: Suppose that $d\geq 4$ is even. Then theorem \ref{explicit} yields that 
\begin{equation} \label{size11}
|S_t| \geq q^{d-1} - q^{\frac{d-2}{2}}.
\end{equation} 
From Theorem \ref{ExplicitGauss}, note that $ \left(G_1(\psi, \chi)\right)^d$ is a real number if $d$ is even. 
Therefore the following two values take the different signs:
\begin{equation}\label{twoterms} - q^{-2}(q-2) \left(G_1(\psi, \chi)\right)^d \quad \mbox{and} \quad q^{-1}\left(G_1(\psi, \chi)\right)^d  
\sum_{r\not=0, 1: t(1-r)^2+r\|x\|=0}1.\end{equation}
Moreover, $$\sum_{r\not=0, 1: t(1-r)^2+r\|x\|=0}1 \leq 2,$$
because the polynomials of degree two have at most two roots. 
Together with this , (\ref{size11}), and (\ref{twoterms}), the second part of Theorem \ref{intersection2} gives 
\begin{align*} |S_t \cap (S_t+x) | \geq& q^{-d}( q^{d-1} - q^{\frac{d-2}{2}})^2-q^{-2}\\
 &- \max \left\{ |q^{-2}(q-2) \left(G_1(\psi, \chi)\right)^d|,\,
|2q^{-1}\left(G_1(\psi, \chi)\right)^d |\right \}\\
=&q^{-d}( q^{d-1} - q^{\frac{d-2}{2}})^2-q^{-2}- \max \left\{ q^{-2}(q-2) q^{\frac{d}{2}},\,2q^{-1}q^{\frac{d}{2}}\right \}\\
=&q^{-d}( q^{d-1} - q^{\frac{d-2}{2}})^2-q^{-2}-2q^{-1}q^{\frac{d}{2}} \\
=&q^{\frac{d-2}{2}} (q^{\frac{d-2}{2}} -2-2q^{-1}), \end{align*}
which is greater than zero if $d\geq 4 , q\geq 3.$ Thus in order to conclude that in dimensions four and higher
the diameter of $G_q^{\Delta}$ is two , it remains to prove Lemma \ref{intersection2}, which shall be done by the following subsection.

\vskip.125in 

\subsection{Proof of Lemma \ref{intersection2}}

For each $x\not=(0,\ldots,0)$ and $t\not=0$, we  have 
\begin{align}\label{convolutionsetup} 
|S_t \cap (S_t+x)|=& \sum_y  S_t(y-x) S_t(y)\nonumber\\
=&\sum_y \sum_m \widehat{S_t}(m) \chi((y-x) \cdot m) S_t(y)\nonumber\\
=&q^{-d}{|S_t|}^2+q^d \sum_{m \not=(0, \dots, 0)} {|\widehat{S_t}(m)|}^2 \chi(-x \cdot m). \end{align}
Now, by (\ref{startaverage}) above, we have 
$$\widehat{S_t}(m)=q^{-d-1} (G_1(\psi,\chi))^d\sum_{s \not=0} \chi \left( \frac{||m||}{-4s}-st \right) \psi^d(s)$$ 
where $\psi$ is the quadratic character of ${\Bbb F}_q^{*}$ with $\psi = \overline{\psi} ,\psi(ab)=\psi(a)\psi(b)$.
Plugging this into the second term in (\ref{convolutionsetup}) above we get 

$$ q^{-2} \sum_{m \not=(0, \dots, 0)} \sum_{u \not=0} \sum_{v \not=0} \psi^d(u)
\overline{\psi^d(v)} \chi \left( \frac{||m||}{-4u}-ut \right) \chi \left( \frac{||m||}{4v}+vt \right) \chi(-x \cdot m)$$
$$=q^{-2} \sum_{m} \sum_{u \not=0} \sum_{v \not=0} \psi^d(uv)
 \chi \left( \frac{||m||}{-4u}-ut \right) \chi \left( \frac{||m||}{4v}+vt \right) \chi(-x \cdot m)$$
$$-q^{-2} \sum_{u \not=0} \sum_{v \not=0}\psi^d(uv) \chi(-t(u-v)) = I +II.$$ 
The second term above is given by
\begin{equation}\label{II} II= -q^{-2} \sum_{u \not=0} \sum_{v \not=0}\psi^d(uv) \chi(-t(u-v)) 
=\left\{\begin{array}{ll} -q^{-1} \quad \mbox{if} \quad d \quad \mbox{is odd}\\
-q^{-2} \quad \mbox{if} \quad d \quad \mbox{is even} \end{array}\right.\end{equation}
This follows from the Gauss sum estimates and the facts that $ \psi^d= \psi $ for $d$ odd, $\psi^d\equiv 1 $ for $d$ even, 
and $ \sum\limits_{s\not=0} \chi(ts) =-1$ for $t\not=0.$
On the other hand, using the changing of variables, $u\to u^{-1}, v\to v^{-1}$, the first term above is written by
$$I= q^{-2} \sum_{m}\sum_{u,v \not=0 : u\not=v} \psi^d(uv) \chi\left(-t\left(\frac{1}{u}-\frac{1}{u}\right)\right) 
\chi\left(\frac{(u-v)}{-4}\|m\|-x\cdot m \right),$$
where we also used the fact that if $u=v$ then the sum in $m \in {\mathbb F}_q^d$ vanishes, because $\sum\limits_{m} \chi(-x\cdot m)=0 $ 
for $x \not=(0,\cdots,0).$ Completing the squares ( see Lemma \ref{complete}) , we have
$$ I = q^{-2} \left(G_1(\psi, \chi)\right)^d \psi^d(-4^{-1})\sum_{u,v \not=0 : u\not=v} \psi^d(uv(u-v))
\chi\left(\frac{t(u-v)}{uv}\right)\chi\left(\frac{\|x\|}{u-v}\right).$$
Note that $ \psi(-4^{-1}) = \psi(-1)$, because $4$ is the square number in ${\mathbb F}_q.$ 
Letting $u=s, vu^{-1}=r$, we see that
$$I= q^{-2} \left(G_1(\psi, \chi)\right)^d \psi^d(-1)\sum_{s\not=0}\sum_{r\not=0,1} 
\psi^d(r(1-r))\psi^d(s) \chi\left(\frac{t(1-r)^2+r\|x\|}{sr(1-r)}\right).$$

Case 1: Suppose that $d$ is even. Then $ \psi^d \equiv 1$. Thus we have
$$I=q^{-2} \left(G_1(\psi, \chi)\right)^d \sum_{s\not=0}\sum_{r\not=0,1}\chi\left(\frac{t(1-r)^2+r\|x\|}{sr(1-r)}\right).$$  
Note that  the sum in $s\not=0$ is $q-1$ if  $t(1-r)^2+r\|x\|=0,$ and $-1$ otherwise. Thus we obtain that
\begin{align*} I=& q^{-2} \left(G_1(\psi, \chi)\right)^d (q-1)\sum_{\substack{r\not=0,1 \\:t(1-r)^2+r\|x\|=0 }}1
-q^{-2} \left(G_1(\psi, \chi)\right)^d\sum_{\substack{r\not=0,1 \\:t(1-r)^2+r\|x\|\not=0 }}1\\
=& q^{-1} \left(G_1(\psi, \chi)\right)^d\sum_{r\not=0,1 :t(1-r)^2+r\|x\|=0 }1-q^{-2}(q-2) \left(G_1(\psi, \chi)\right)^d .\end{align*}
From this estimate , (\ref{convolutionsetup}), and (\ref{II}), the second part of Lemma \ref{intersection2} immediately follows.\\

Case 2: Suppose that $d$ is odd. Then $ \psi^d =\psi$. Since $ \sum\limits_{s\not=0} \psi(s) =0$, we see that 
the sum in $s\not =0$ is zero if $t(1-r)^2+r\|x\|=0.$ Thus we may assume that $t(1-r)^2+r\|x\|\not=0$. Using  the changing of variables,
$s^{-1}r^{-1} (1-r)^{-1} ( t(1-r)^2+r\|x\|) \to s $ and the facts that $\psi(r^2)=1, \psi(r)=\psi(r^{-1})$, we see that
$$ I =q^{-2} \left(G_1(\psi, \chi)\right)^{d+1} \psi(-1)\sum_{ r\not=0,1: t(1-r)^2+r\|x\|\not=0 }\psi\left(t(1-r)^2+r\|x\|\right).$$
Using this estimate together with (\ref{convolutionsetup}) and (\ref{II}), the first part of Lemma \ref{intersection2}  is proved and so the proof is complete.

\vskip.125in 
\subsection{The diameter of $G_q^{\Delta}$ in two dimension is never two}\label{nevertwo} As we claim in the statement of the theorem, 
we can show that for any field ${\Bbb F}_q$ with $q\not= 3$, the diameter of $G_q^{\Delta}$ in dimension two is indeed  three and never two . 
To prove this it suffices to show that for each $t\not=0 $ there exists $x \in {\Bbb F}_q^2$ 
such that the circle $S_t$ and its translate by $x$ do not intersect.\\

Case 1. Suppose that $q$ is any power of odd prime $p\equiv 1\, (mod\,4)$, 
or $q$ is even power of odd prime $p\equiv 3\, (mod\,4).$ Then  $q\equiv 1\, (mod\,4)$ which says that $-1$ is a square number in ${\mathbb F}_q$ so that
$\psi(-1) =1.$ By this and Theorem \ref{explicit}, if $d=2$ then $|S_t|= q-1.$ Moreover we see from Theorem \ref{ExplicitGauss} that $ \left(G_1(\psi,\chi)\right)^2=q.$ 
Using the second part of Lemma \ref{intersection2}, we therefore obtain that for each $x\not= (0,0), t\not=0,$
\begin{align}\label{zero} |S_t \cap (S_t + x)|=&q^{-2}(q-1)^2-q^{-2}-q^{-1}(q-2) + \sum_{r\not= 0,1: t(1-r)^2 +r\|x\|=0} 1\nonumber\\
                                  =& \sum_{r\not= 0,1: t(1-r)^2 +r\|x\|=0} 1.\end{align}
If we choose $ x=(1,i)\in {\mathbb F}_q^2$ with $ i^2 = -1$ then $ \|x\|=0$ and the sum in (\ref{zero}) vanishes.
Thus two circles $S_t$ and $(S_t +x)$ are disjoint.\\

Case 2: Suppose that $q$ is an odd power of odd prime $p\equiv 3 \,(mod\,4).$ Then $q\equiv 3\, (mod\,4)$ and so $\psi(-1)=-1$, because $-1$ is not a square number in
${\mathbb F}_q.$ Together with this, Theorem \ref{explicit} implies that $|S_t|= q+1$ if $d=2.$ In addition, we see from Theorem \ref{ExplicitGauss} that 
$ \left(G_1(\psi,\chi)\right)^2=-q.$ Thus the second part of Lemma \ref{intersection2} yields that for each $x\not= (0,0), t\not=0,$
\begin{align*}|S_t \cap (S_t + x)|=&q^{-2}(q+1)^2-q^{-2}+q^{-1}(q-2) - \sum_{r\not= 0,1: t(1-r)^2 +r\|x\|=0} 1\\
                                 =&2-\sum_{r\not= 0,1: t(1-r)^2 +r\|x\|=0} 1.\end{align*}
It therefore suffices to show that for each $t\not=0,$ there exists $ x\not =(0,0)$ such that 
\begin{equation}\label{aim} D_t(x)=\sum_{r\not= 0,1: t(1-r)^2 +r\|x\|=0} 1 =2.\end{equation}
Observe that $\|x\|\not=0$ if $x\not =(0,0)$, because we have assumed that $-1$ is not a square number.
Thus if $x\not =(0,0),$ then $t(1-r)^2 +r\|x\|\not=0$ for $r=0, 1,$ because $t\not=0.$ From this observation, we see that for each $ x\not =(0,0),$
\begin{align*} D_t(x)=&\sum_{r\in {\mathbb F}_q: t(1-r)^2 +r\|x\|=0} 1\\
                     =& q^{-1} \sum_{s,r \in {\mathbb F}_q} \chi\left(s (t(1-r)^2 +r\|x\|)\right)\\
                     =& 1+ q^{-1} \sum_{r\in {\mathbb F}_q} \sum_{s\not=0} \chi\left(s (t(1-r)^2 +r\|x\|)\right)\\
                     =& 1+q^{-1}\sum_{s\not=0}\sum_{r\in {\mathbb F}_q}\chi\left(str^2 +(s\|x\|-2st)r\right) \chi(st)\\
                     =& 1+ q^{-1}G_1(\psi,\chi)\sum_{s\not=0}\psi(st) \chi\left(\frac{s(\|x\|-2t)^2}{-4t}\right) \chi(st), \end{align*}
where the last equality can be obtained by the completing square methods in Lemma \ref{complete}.
Using the changing of variables, $s/(-4t) \to s$, we have 
$$ D_t(x)= 1+q^{-1}G_1(\psi,\chi)\psi(-1) \sum_{s\not=0}\psi(s)\chi\left( s \left((\|x\|-2t)^2-4t^2\right)\right),$$
because $ \psi $ is the quadratic multiplicative character of ${\mathbb F}_q^*$ and so $\psi(4s^2)=1.$
Here, recall that we have assumed that $-1$ is not a square number and so $\psi(-1)=-1 , \|x\|\not=0$ for $x\not=(0,0).$
In addition, we assume that $ \|x\|\not=4t.$ Then  $(\|x\|-2t)^2-4t^2$ can not be zero. Thus we apply the changing of variables,
$s \left((\|x\|-2t)^2-4t^2\right)\to s$ and we obtain that
$$ D_t(x)=1-q^{-1}\left(G_1(\psi,\chi)\right)^2 \psi\left((\|x\|-2t)^2-4t^2)\right),$$
where we used the fact that $ \psi(s)=\psi(s^{-1})$ for $s\not=0.$ 
By Theorem \ref{ExplicitGauss} and our assumption in Case 2, observe that $  \left(G_1(\psi,\chi)\right)^2 =-q.$
Thus if $\|x\|\not=0, 4t$, then $D_t(x)$ above takes the following form.
$$ D_t(x)=1+\psi\left(\|x\|^2-4t\|x\|\right).$$
In order to prove (\ref{aim}), it therefore suffices to show that 
for each $ t\not=0,$ there exists $ x\in {\mathbb F}_q^2 $ with $\|x\|\not=0, 4t$ such that 
$$  \psi\left(\|x\|^2-4t\|x\|\right) =1.$$
By contradiction, we assume that for all $x\in {\mathbb F}_q^2 $ with $\|x\|\not=0, 4t$, 
\begin{equation}\label{realfalse}\psi\left(\|x\|^2-4t\|x\|\right) \not=1.\end{equation}
Since $\|x\|^2-4t\|x\|\not=0 $ for $ x\in {\mathbb F}_q^2 $ with $\|x\|\not=0, 4t,$ and $\psi$ is the quadratic character of ${\mathbb F}_q^*$, we see that
$\psi(s)$ for $s\not=0$ takes $+1$ or $-1.$ Moreover, observe from Theorem \ref{explicit} that for each $s\not=0,$ there exists $x\in {\mathbb F}_q^2$
such that $ \|x\|=s.$ Thus (\ref{realfalse}) implies that 
$$ \sum_{\|x\| \in {\mathbb F}_q\setminus \{0, 4t\}}\psi\left(\|x\|^2-4t\|x\|\right) = -(q-2).$$ 
By defining $\psi(0)=0$ we extend the quadratic character $\psi$ of ${\mathbb F}_q^* $ to the quadratic character of ${\mathbb F}_q.$
Then we have  
\begin{equation}\label{anotherfalse} \sum_{\|x\| \in {\mathbb F}_q}\psi\left(\|x\|^2-4t\|x\|\right) = -(q-2).\end{equation}
However, this is impossible if $q\geq 5$ due to the following theorem ( See \cite{LN97}, P.225 ).
\begin{theorem}\label{goodthm}
Let $\psi$ be a multiplicative character of ${\mathbb F}_q$ of order $k>1$ and let $g\in {\mathbb F}_q[x]$ be a monic polynomial of positive degree that is
not an $k-$th power of a polynomial. Let $e$ be the number of distinct roots of $g$ in its splitting field over ${\mathbb F}_q.$ Then for every $s\in {\mathbb F}_q$
we have
$$\left| \sum_{t\in {\mathbb F}_q} \psi(sg(t)) \right|\leq (e-1) q^{1/2}.$$
\end{theorem}
To see that $(\ref{anotherfalse})$ is false if $q\geq 5$, note from Theorem \ref{goodthm} that 
 $$\left|\sum_{\|x\| \in {\mathbb F}_q}\psi\left(\|x\|^2-4t\|x\|\right)\right|\leq q^{\frac{1}{2}}.$$
Thus the proof is complete.

\subsection{The two-dimensional case and the three-dimensional case}\label{lastproof} 
In this subsection, we prove that the diameter of $G_q^{\Delta}$ in two dimension is three unless $q=3,5,9, 13,$ and 
the diameter of $G_q^{\Delta}$ in three dimension is less than equal to three.
For each $a,b \in {\mathbb F}_q^d$ with $a\not= b$, it suffices to show that
$$ |\{(x,y)\in (S_r+a) \times (S_r+b): \|x-y\| =r\not=0 \}| >0.$$
We have
\begin{align}\label{intersection3}&|\{(x,y)\in (S_r+a) \times (S_r+b): \|x-y\| =r\not=0 \}|\nonumber \\
=& \sum_{x,y} S_r(x)S_r(y)S_r(x-y+a-b)\nonumber \\
=& \sum_{x,y}\sum_{m} S_r(x)S_r(y) \chi\left((x-y+a-b)\cdot m\right) \widehat{S_r}(m)\nonumber \\
=&q^{-d}|S_r|^3 + q^{2d} \sum_{m\not=(0,\ldots,0)} \chi((a-b)\cdot m) \widehat{S_r}(m) |\widehat{S_r}(m)|^2.\end{align}
We shall estimate the second term above. By (\ref{startaverage}), recall that $\widehat{S_r}(m)$  is given by
$$\widehat{S_r}(m)=q^{-d-1} (G_1(\psi,\chi))^d\sum_{s \not=0} \chi \left( \frac{||m||}{-4s}-sr \right) \psi^d(s).$$
Plugging this into the second part in (\ref{intersection3}), we have
$$ q^{-3} \left(G_1(\psi, \chi)\right)^d\sum_{m\not=(0,\ldots,0)}\sum_{u,v,w\not=0}\psi^d(uvw)
\chi((a-b)\cdot m)\chi\left( \frac{\|m\|}{-4}\left(\frac{1}{u}-\frac{1}{v}+\frac{1}{w}\right)\right)\chi(-r(u-v+w))$$
$$=q^{-3} \left(G_1(\psi, \chi)\right)^d\sum_{m}\sum_{u,v,w\not=0}\psi^d(uvw)
\chi((a-b)\cdot m)\chi\left( \frac{\|m\|}{-4}\left(\frac{1}{u}-\frac{1}{v}+\frac{1}{w}\right)\right)\chi(-r(u-v+w))$$
$$- q^{-3} \left(G_1(\psi, \chi)\right)^d\sum_{u,v,w\not=0}\psi^d(uvw)\chi(-r(u-v+w)) =I +II.$$
The second term $II$ above is given by 
\begin{equation}\label{II-1} \left\{\begin{array}{ll} q^{-3} \left(G_1(\psi, \chi)\right)^d \quad &\mbox{if} \quad d \quad \mbox{is even}\\
                                  -q^{-3} \left(G_1(\psi, \chi)\right)^{d+3} \psi(r)   \quad & \mbox{if} \quad d \quad \mbox{is odd} \end{array}\right.\end{equation}
 This easily follows from  properties of the quadratic character $\psi$, definition of the Gauss sum $G_1(\psi, \chi)$, and $ \sum\limits_{s\not=0} \chi(rs)=-1$ for $r\not=0.$
Let us estimate the first term $I$ above. Using the changing of variables, $u^{-1}\to u, v^{-1}\to v, w^{-1}\to w,$ the first term $I$ above is given by
$$ q^{-3} \left(G_1(\psi, \chi)\right)^d\sum_{m}\sum_{u,v,w\not=0}\psi^d(uvw)
\chi((a-b)\cdot m)\chi\left( \frac{\|m\|}{-4}(u-v+w)\right)\chi\left(-r\left(\frac{1}{u}-\frac{1}{v}+\frac{1}{w}\right)\right).$$  
Since $a\not=b$, the sum in $m \in {\mathbb F}_q^d$ vanishes if $ u-v+w =0.$ Thus we may assume that $ u-v+w \not=0$.
Therefore using  Lemma \ref{complete}, the first term $I$ above takes the form
\begin{equation}\label{simpleI} q^{-3} \left(G_1(\psi, \chi)\right)^{2d} \sum_{\substack{u,v,w\not=0 \\: u-v+w\not=0}} 
\psi^d(uvw)\psi^d(-(u-v+w)) \chi\left(\frac{\|a-b\|}{u-v+w} \right)\chi\left(-r\left(\frac{1}{u}-\frac{1}{v}+\frac{1}{w}\right)\right) ,\end{equation}   
where we used $ \psi(4^{-1})=1$ , because $4$ is the square number and $\psi$ is the quadratic character of ${\mathbb F}_q^*.$                   
Now we estimate the term $I$ above in the cases when $d=3$ and $d=2$.\\

Case A: The dimension $d$ is three. Then the term $I$ is dominated by 
$$ |I|\leq \sum_{u,v,w \not=0: u-v+w \not=0} 1,$$
where we used the fact that the magnitude of the Gauss sum $G_1(\psi, \chi)$ is exactly $q^{\frac{1}{2}}.$
We claim that 
\begin{equation}\label{claim}\sum_{u,v,w \not=0: u-v+w \not=0} 1 = (q-1)^2 + (q-1)(q-2)^2.\end{equation}
The claim follows from the following observation:
if we fix $u\not=0$ which has $d-1$ choices, then we may choose $v\not=0$ such that either $u=v$ or $u\not=v.$
In case $u=v$, $v$ has only one choice which depends on the choice of $u$ and then we can choose $w\not=0 $ which has $q-1$ choices with $u-v+w\not=0.$
One the other hand, if we choose $v\not=0 $ with $u\not=v$, which has $q-2$ choices, then we have $q-2$ choices for $w\not=0$ so that $u-v+w\not=0.$
Thus the claim holds. From (\ref{intersection3}), (\ref{II-1}), and (\ref{claim}) above, we obtain that if $d=3$ then
$$ |\{(x,y)\in (S_r+a) \times (S_r+b): \|x-y\| =r\not=0 \}| $$ 
$$\geq q^{-3}|S_r|^3 -|q^{-3} \left(G_1(\psi, \chi)\right)^6 \psi(r) |-\left((q-1)^2 + (q-1)(q-2)^2\right)$$
$$\geq q^{-3} (q^2-q)^3 -1 -(q-1)^2 - (q-1)(q-2)^2 = (q-1)\left(q-2- (q-1)^{-1}\right),$$
which is greater than zero if $q\geq 3.$ This proves that the diameter of $G_q^{\Delta}$ in three dimension is less than equal to three.
\begin{remark} In three dimension , the diameter of $G_q^{\Delta}$ depends on both the finite field ${\mathbb F}_q$ and 
the choice of the radius $r\not=0$ of $S_r.$ In fact, by estimating the sum in the first part of Lemma \ref{intersection2}, we can show that
the diameter of  $G_q^{\Delta}$ in three dimension is two if $ \psi(-r)=1$, and three otherwise. This can be done by the similar arguments as in Subsection \ref{nevertwo}
and this subsection.
\end{remark}

Case B: The dimension $d$ is two. Then the term $I$ in (\ref{simpleI}) above takes the form
$$ I= q^{-3} \left(G_1(\psi, \chi)\right)^4 \sum_{\substack{u,v,w\not=0 \\: u-v+w\not=0}} \chi\left((u-v+w)^{-1} \|a-b\|\right)\chi\left(-r(u^{-1}-v^{-1}+w^{-1})\right) .$$
Fix $u\not=0.$ Putting $ u^{-1}v=s, u^{-1}w=t$, we see that 
$$ I = q^{-3} \left(G_1(\psi, \chi)\right)^4 \sum_{u\not=0}\sum_{s\not=0}\sum_{\substack{t\not=0: s-t\not=1}} 
\chi\left(-r\left(\frac{1}{u}+ \frac{s-t}{ust}\right)\right) \chi\left(\frac{\|a-b\|}{u(1-s+t)}\right).$$
Using the changing of variables, $u^{-1}\to u$, we have
$$I= q^{-3} \left(G_1(\psi, \chi)\right)^4 \sum_{u\not=0}\sum_{s\not=0}\sum_{\substack{t\not=0: s-t\not=1}} 
 \chi\left(\left(-r+\frac{-rs+rt}{st}+ \frac{\|a-b\|}{1-s+t}\right) u\right) .$$   
Note that the sum in $u\not=0$ is $-1$ if $-r+(-rs+rt)/st + \|a-b\|/(1-s+t)\not=0,$ and $q-1$ otherwise. Thus the term $I$ can be written by
\begin{align*} I =& -q^{-3} \left(G_1(\psi, \chi)\right)^4\sum_{s\not=0}\sum_{\substack{t\not=0: s-t\not=1\\ -r+(-rs+rt)/st + \|a-b\|/(1-s+t)\not=0}} 1\\
                  &+ q^{-3} \left(G_1(\psi, \chi)\right)^4 (q-1) \sum_{s\not=0}\sum_{\substack{t\not=0: s-t\not=1\\ -r+(-rs+rt)/st + \|a-b\|/(1-s+t)=0}} 1\\
                 =& q^{-2}\left(G_1(\psi, \chi)\right)^4\sum_{s\not=0}\sum_{\substack{t\not=0: s-t\not=1\\ -r+(-rs+rt)/st + \|a-b\|/(1-s+t)=0}} 1\\
                  &- q^{-3}\left(G_1(\psi, \chi)\right)^4\sum_{s\not=0}\sum_{t\not=0: s-t\not=1} 1.\end{align*}
We now claim that 
$$\sum_{s\not=0}\sum_{t\not=0: s-t\not=1} 1= (q-2)^2 +(q-1).$$
To see this, we write the term above into two parts as follows.
$$ \sum_{s\not=0}\sum_{t\not=0: s-t\not=1} 1 
= \sum_{s\not=0,1}\sum_{t\not=0: s-t\not=1} 1 + \sum_{t\not=0: 1-t\not=1}1.$$
Then it is clear that
$$\sum_{t\not=0: 1-t\not=1}1=q-1.$$
On the other hand, we see that 
$$\sum_{s\not=0,1}\sum_{t\not=0: s-t\not=1} 1 =(q-1)^2,$$
because whenever we fix $s\not=0, 1$ which has $q-2$ choices, we have $q-2 $ choices of $t\not=0$ with $s-t\not=1.$
By this, the claim is complete. Thus the term $I$ is given by
$$ I=q^{-2}\left(G_1(\psi, \chi)\right)^4\sum_{s\not=0}\sum_{\substack{t\not=0: s-t\not=1\\ -r+(-rs+rt)/st + \|a-b\|/(1-s+t)=0}} 1$$
$$- q^{-3}\left((q-2)^2 +(q-1)\right) \left(G_1(\psi, \chi)\right)^4.$$
From this, (\ref{intersection3}), and (\ref{II-1}), we obtain that if $d=2$ then
$$ |\{(x,y)\in (S_r+a) \times (S_r+b): \|x-y\| =r\not=0 \}| $$ 
$$= q^{-2}|S_r|^3 + q^{-3} \left(G_1(\psi, \chi)\right)^2 - q^{-3}\left((q-2)^2 +(q-1)\right) \left(G_1(\psi, \chi)\right)^4$$
$$+q^{-2}\left(G_1(\psi, \chi)\right)^4\sum_{s\not=0}\sum_{\substack{t\not=0: s-t\not=1\\ -r+(-rs+rt)/st + \|a-b\|/(1-s+t)=0}} 1.$$
By Theorem \ref{explicit}) and Theorem \ref{ExplicitGauss}, we see that $|S_r|=q-\psi(-1)$ if $d=2$, and $G^4(\psi,\chi)=q^2$ respectively.
Thus we aim to show that the following value is positive.
\begin{equation}\label{formulathree} |\{(x,y)\in (S_r+a) \times (S_r+b): \|x-y\| =r\not=0 \}| \end{equation} 
$$=q^{-2}(q-\psi(-1))^3 + q^{-3} \left(G_1(\psi, \chi)\right)^2 - q^{-1}\left(q^2-3q+3\right) $$
$$+\sum_{s\not=0}\sum_{\substack{t\not=0: s-t\not=1\\ -r+(-rs+rt)/st + \|a-b\|/(1-s+t)=0}} 1.$$

Case B-I: Suppose that $q=p^l$ for some odd prime $p\equiv 3 \, (mod\,4)$ with $l$ odd.
Then $q\equiv 3\,(mod\,4)$ which means that $-1$ is not a square number in ${\mathbb F}_q$ so that $ \psi(-1)=-1.$
We also note from Theorem \ref{ExplicitGauss} that $ G^2(\psi, \chi) =-q.$ Thus the term in (\ref{formulathree}) can be estimated as follows.
$$|\{(x,y)\in (S_r+a) \times (S_r+b): \|x-y\| =r\not=0 \}|$$
$$\geq q^{-2}(q+1)^3-q^{-2}- q^{-1}\left(q^2-3q+3\right) =6$$
which is greater than zero as wanted.\\

Case B-II: Suppose that  $q=p^l$ for some odd prime $p\equiv 3\, (mod\,4)$ with $l$ even, or  $q=p^l$ with $p\equiv 1\,(mod\,4).$
Then $q\equiv 1\,(mod\,4)$ which implies that $-1$ is a square number in ${\mathbb F}_q$ so that $\psi(-1)=1.$
Moreover $ G^2(\psi,\chi) =q $ by Theorem \ref{ExplicitGauss}. From these observations, the term in (\ref{formulathree}) is given by
$$|\{(x,y)\in (S_r+a) \times (S_r+b): \|x-y\| =r\not=0 \}|$$
$$ =q^{-2}(q-1)^3 +q^{-2}-q^{-1}(q^2-3q+3) + R(a,b,r) = R(a,b,r)$$
where $$R(a,b,r)= \sum_{s\not=0}\sum_{t\not=0: s-t\not=1 , T(s,t,a,b,r)=0} 1$$
with 
$$T(s,t,a,b,r)=-r+(-rs+rt)/st + \|a-b\|/(1-s+t).$$
To complete the proof, it suffices to show that $ R(a,b,r) >0.$\\

Case B-II-1: Suppose that $\|a-b\|=0.$ Then we have 
$$ R(a,b,r) = \sum_{s\not=0}\sum_{\substack{t\not=0: s-t\not=1,\\ -st-s+t=0}}1.$$
If $q\not=3^l$ for $l$ even $(Char\,{\mathbb F}_q\not=3)$, then it is clear that 
$ R(a,b,r)\geq 1$, because if we choose $s=-1, t=-2^{-1}$ then $s-t\not\equiv 1$ and $-st-s+t \equiv 0.$
Thus we may assume that $q=3^l$ with $l$ even.
Since each finite field ${\mathbb F}_{3^2}$ can be considered as a subfield of any finite field ${\mathbb F}_{3^l}$ with $l$ even up to isomorphism,
it is enough to show that $ R(a,b,r) \geq 1 $ for a fixed finite field ${\mathbb F}_{3^2}$ with $ 9$ elements.
Consider the following finite field ${\mathbb F}_{3^2}$ with $9$ elements.
$${\mathbb F}_{3^2} \cong {\mathbb Z}_3 [i]/ (i^2+1) \cong \{\alpha+\beta i: \alpha, \beta \in {\mathbb Z}_3\},$$ 
where $i^2 =-1.$ 
Taking $s=i, t= \frac{i-1}{2}$, we see that $s-t \not\equiv 1$ and $-st-s+t\equiv 0.$ Thus we conclude that $R(a,b,r)\geq 1$ as desired.\\

Case B-II-2. Assume that $\|a-b\|\not=0.$ Letting $ c=\frac{\|a-b\|}{r} \not=0$, we have
$$ R(a,b,r)=\sum_{s\not=0}\sum_{t\not=0: s-t\not=1 , T^*(s,t,c)=0} 1$$
where $T^*(s,t,c)$ is defined by
\begin{align}\label{defT*} T^*(s,t,c) =& (c-3)st+s^2t-st^2-s+t+s^2+t^2\nonumber\\
                          =& (t+1)s^2+ \left(t(c-3)-t^2-1\right)s +t+t^2.\end{align}
Splitting $R(a,b,r)$ into two parts as below and using the simple properties of summation notation, $R(a,b,r)$ takes the following forms.
\begin{align*}R(a,b,r)=&\sum_{s\not=0,1}\sum_{t\not=0: s-t\not=1 , T^*(s,t,c)=0} 1 + \sum_{t\not=0 : c-1=0}1\\
=&\sum_{s\not=0,1}\sum_{t\not=0:  T^*(s,t,c)=0} 1 -\sum_{s\not=0,1 : c=0}1 + \sum_{t\not=0 : c-1=0}1 \\
=&\sum_{s\not=0}\sum_{t\not=0:  T^*(s,t,c)=0} 1-\sum_{t\not=0 : c-1=0}1-\sum_{s\not=0,1 : c=0}1+ \sum_{t\not=0 : c-1=0}1\\
=&\sum_{s\not=0}\sum_{t\not=0:  T^*(s,t,c)=0} 1,\end{align*}
where we used  $\sum\limits_{s\not=0,1 : c=0}1=0,$ because $c\not=0.$ We have
\begin{align}\label{split}R(a,b,r)=& q^{-1} \sum_{s,t\not=0} \sum_{k} \chi\left( k T^*(s,t,c)\right)\nonumber\\
                      =& q^{-1} \sum_{s,t\not=0} \sum_{k\not=0} \chi\left( k T^*(s,t,c)\right) +q^{-1}(q-1)^2 \nonumber\\
                      =& q^{-1} \sum_{t,k\not=0} \sum_{s\in {\mathbb F}_q} \chi\left( k T^*(s,t,c)\right)
                        -q^{-1}\sum_{t,k\not=0} \chi(tk+t^2k) +q^{-1}(q-1)^2 \nonumber\\
                      =& q^{-1} \sum_{t,k\not=0} \sum_{s\in {\mathbb F}_q} \chi\left( k T^*(s,t,c)\right)
                        -q^{-1}+q^{-1}(q-1)^2, \end{align}
where the last equality follows from the following observation.
\begin{align*}\sum_{t,k\not=0} \chi(tk+t^2k) =& \sum_{t\not=0,-1}\sum_{k\not=0} \chi(t(t+1)k) + \sum_{k\not=0}1\\
                                             =& -(q-2)+(q-1) =1.\end{align*}
Splitting the sum in (\ref{split}) into two parts as below, we obtain that
$$ R(a,b,r) =q^{-1}\sum_{t\not=0,-1}\sum_{k\not=0} \sum_{s\in {\mathbb F}_q} \chi\left( k T^*(s,t,c)\right)
+q^{-1}\sum_{k\not=0} \sum_{s\in {\mathbb F}_q} \chi\left( (1-c)ks\right)-q^{-1}+q^{-1}(q-1)^2.$$
By the orthogonality relations for non-trivial additive character $\chi$, the second term above is given by
$$q^{-1}\sum_{k\not=0} \sum_{s\in {\mathbb F}_q} \chi\left( (1-c)ks\right) = (q-1) \,\delta_0(1-c) \geq 0,$$
where $\delta_0(u)=1$ if $u=0$, and $0$ otherwise. In order to estimate the first term above, 
 recall from (\ref{defT*}) that 
$$   k T^*(s,t,c) =k(t+1)s^2+ k\left(t(c-3)-t^2-1\right)s +kt+kt^2$$
and then apply the complete square methods ( see Lemma \ref{complete} ). It follows that
$$ R(a,b,r)\geq q^{-1}G_1(\psi,\chi)\sum_{t\not=0,-1}\sum_{k\not=0} \psi\left((t+1)k\right) \chi\left(\frac{k\left((c-3)t-t^2-1\right)^2}{-4(t+1)} \right)\chi(t(t+1)k)
+q-2.$$
Using the changing of variables, $\frac{k}{4(t+1)} \to k$ and the fact that $ \psi\left(4(t+1)^2\right) =1$, we see that
$$R(a,b,r)\geq q^{-1}G_1(\psi,\chi)\sum_{t\not=0,-1}\sum_{k\not=0} \psi(k) \chi( g(t,c) k) + q-2 $$
where $g(t,c)$ is given by
$$ g(t,c) = 4t(t+1)^2-\left((c-3)t-t^2-1\right)^2.$$
Note that the sum in $k\not=0$ is zero if $ g(t,c)=0$. Thus we may assume that $g(t,c)\not=0.$ Thus using the changing variables,
$g(t,c) k \to k$, we see that
\begin{equation}\label{almostfinal}R(a,b,r)\geq q^{-1}\left(G_1(\psi,\chi)\right)^2\sum_{t\not=0,-1 : g(t,c)\not=0} \psi(g(t,c)) +q-2,\end{equation}
where we used that $\psi(s) = \psi(s^{-1})$ for $s\not=0.$
Here, recall that we has assumed that $q=p^l$ for some odd prime $p\equiv 3\, (mod\,4)$ with $l$ even, or  $q=p^l$ with $p\equiv 1\,(mod\,4).$
By Theorem \ref{ExplicitGauss}, we therefore see that $\left(G_1(\psi,\chi)\right)^2= q.$ From this and (\ref{almostfinal}), we see that
$R(a,b,r)=0$ only if $ \sum\limits_{t\not=0,-1 : g(t,c)\not=0} \psi(g(t,c))= -(q-2).$ Since $\psi$ is the quadratic character of ${\mathbb F}_q^*$, 
the number $\psi(g(t,c))$ takes $+1$ or $-1$. Thus if $\sum\limits_{t\not=0,-1 : g(t,c)\not=0} \psi(g(t,c))= -(q-2)$ happens then it must be true that
$\psi(g(t,c))=-1$ for all $t\not=0,-1.$ This implies that $g(t,c)$  is not a square number for all $t\not=0,-1,$ and the following estimate holds
\begin{equation}\label{false} \left|\sum_{t\in {\mathbb F}_q} \psi(g(t,c))\right|\geq |-(q-2)|-2 =q-4.\end{equation}
However, by Theorem \ref{goodthm}, it must be true that
$$ \left|\sum_{t\in {\mathbb F}_q} \psi(g(t,c))\right|\leq 3 q^{\frac{1}{2}}. $$
because $g(t,c)$ is the polynomial of degree four in terms of $t$ variables.
Thus if $q\geq 17$ then the inequality in (\ref{false}) is not true and so we conclude that  under the assumptions in Case $B-II$, if $q\geq 17$ then
$$|\{(x,y)\in (S_r+a) \times (S_r+b): \|x-y\| =r\not=0 \}| = R(a,b,r) >0.$$
Combining this and results from Case $B-I$, we finish proving that the diameter of $G_q^{\Delta}$ in two dimension is three if $q\not=3,5,9,13,$
because  the diameter of $G_q^{\Delta}$ in two dimension is never two if $q\not=3.$

\end{document}